\documentclass[12pt]{article}
\usepackage{amsmath}
\usepackage{amssymb}
\usepackage{latexcad}
\usepackage{url}
\usepackage{color}
\definecolor{goodgreen}{rgb}{0,.6,0}
\definecolor{goodred}{rgb}{.8,0,0}
\definecolor{goodblue}{rgb}{0,0,.7}
\usepackage[dvips,colorlinks,citecolor=goodgreen,linkcolor=goodred,urlcolor=goodblue,bookmarks]{hyperref}
\usepackage[right,displaymath,mathlines]{lineno} 
\usepackage{newproof}

\def\Ff{\mathsf F}                                   
\def\FF{\mathop{\cal F}\nolimits}                    
\def\F{\Phi}                                         
\def\La{\Lambda}                                     
\def\e{\varepsilon}                                  
\def\beq{\begin{equation}}                           
\def\eeq{\end{equation}}                             
\def\beqq{\begin{eqnarray}}                          
\def\eeqq{\end{eqnarray}}                            
\def\liml{\mathop{\lim}\limits}                      
\def\suml{\mathop{\sum}\limits}                      
\def\adj{\mathop{\rm adj}\nolimits}                  
\def\Pr{\mathop{\rm Pr\,}\nolimits}                  
\newcommand{\abss}[1]{\left|#1\right|}               
\newcommand{\abs}[1]{|\!#1\!|}                       
\def\cdc{,\ldots,}                                   
\def\on{\overline{1,n}}                              
\newtheorem{thm}{Theorem}{\bfseries}{\itshape}       
\newtheorem{prop}{Proposition}{\bfseries}{\itshape}  
\newtheorem{lemma}{Lemma}{\bfseries}{\itshape}       
\newtheorem{corol}{Corollary}{\bfseries}{\itshape}   
\title{SPANNING FORESTS AND THE GOLDEN RATIO\footnote{To appear in Disc. Appl. Math. (2007),
\url{http://dx.doi.org/10.1016/j.dam.2007.08.030}}}

\author{Pavel Chebotarev\footnote{E-mail addresses: {\tt chv@member.ams.org, upi@ipu.ru
}}\\
{\normalsize Trapeznikov Institute of Control Sciences of the Russian Academy of Sciences}\\
{\normalsize 65 Profsoyuznaya Street, Moscow 117997, Russia}
}
\date{}

\begin{document}

\maketitle
\unitlength 1.50mm 

\begin{abstract}
For a graph $G$, let $f_{ij}$ be the number of spanning rooted forests in which vertex $j$ belongs to a tree
rooted at~$i$. In this paper, we show that for a path, the $f_{ij}$'s can be expressed as the products of
Fibonacci numbers; for a cycle, they are products of Fibonacci and Lucas numbers. The {\em doubly stochastic
graph matrix\/} is the matrix $F=\frac{(f_{ij})_{n\times n}}{f}$, where $f$ is the total number of spanning
rooted forests of $G$ and $n$ is the number of vertices in~$G$. $F$~provides a proximity measure for graph
vertices. By the matrix forest theorem, $F^{-1}=I+L$, where $L$ is the Laplacian matrix of~$G$. We show that
for the paths and the so-called T-caterpillars, some diagonal entries of $F$ (which provide a measure of the
self-connectivity of vertices) converge to $\phi^{-1}$ or to $1-\phi^{-1}$, where $\phi$ is the golden ratio,
as the number of vertices goes to infinity. Thereby, in the asymptotic, the corresponding vertices can be
metaphorically considered as ``golden introverts'' and ``golden extroverts,'' respectively. This metaphor is
reinforced by a Markov chain interpretation of the doubly stochastic graph matrix, according to which $F$
equals the overall transition matrix of a random walk with a random number of steps on~$G$.
\bigskip

\bigskip


\end{abstract}

\section{Introduction}

Let $G=(V,E)$ be a simple graph with vertex set $V=V(G)$, $|V|=n$, and edge set $E=E(G)$.
Suppose that $n\ge2$.

A {\em spanning rooted forest of\/} $G$ is any spanning acyclic subgraph of $G$ with a single vertex
(a~{\em root}) marked in each tree.

Let $f_{ij}=f_{ij}(G)$ be the number of spanning rooted forests of $G$ in which vertices $i$ and $j$ belong
to the same tree rooted at~$i$. The matrix $(f_{ij})_{n\times n}$ is the {\em matrix of spanning rooted
forests of\/}~$G$. Let $f=f(G)$ be the total number of spanning rooted forests of~$G$.

The matrix $F=\frac{(f_{ij})_{n\!\times\! n}}{f}$ is referred to as the {\em doubly stochastic graph
matrix\/}~\cite{Merris97,Merris98,XDZhangWu05,XDZhang05} or the {\em matrix of relative connectivity via forests}.
By the matrix forest theorem~\cite{CheSha95,CheSha97,CheAga02ap},
\beq
\label{mft1}
F^{-1}=I+L
\eeq
and
\beq
\label{mft2}
f=\det(I+L),
\eeq
where $L$ is the {\em Laplacian matrix\/} of $G$, i.e. $L=D-A$, $A$ being the adjacency matrix of $G$ and $D$
the diagonal matrix of vertex degrees of~$G$. Most likely, the matrix $(I+L)^{-1}=F$ was first considered
in~\cite{GolenderDrboglav81}. Chaiken \cite{Chaiken82} used the matrix $\adj(I+L)=(f_{ij})_{n\times n}$
for coordinatizing linking systems of strict gammoids. The $(i,j)$ entry of $F$ can be considered as a measure of
proximity between vertices $i$ and $j$ in $G$; the $(i,i)$ entry measures the self-connectivity of vertex~$i$.

A {\em path\/} is a connected graph in which two vertices have degree~1 and the remaining vertices have
degree~2. Let $P_n$ be the path with $V(P_n)=\{1,2,\ldots,n\}$ and
$E(P_n)=\{(1,2),(2,3),\ldots,(n-1,n)\}$, see Fig.~\ref{fig1}(a).
\unitlength 1.12mm  
\vspace{-1em}
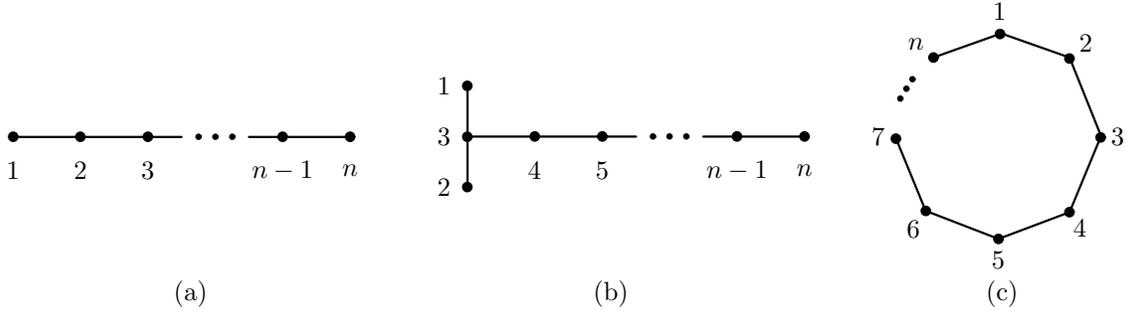
\begin{figure}[htb] 
\hspace{-3em}
\begin{picture}(157.8,46)
\thicklines
\drawthickdot{19.15}{25.04}
\drawthickdot{27.15}{25.04}
\drawthickdot{35.15}{25.04}
\drawthickdot{51.15}{25.04}
\drawthickdot{59.15}{25.04}
\drawpath{19.15}{25.04}{35.15}{25.04}
\drawpath{35.15}{25.04}{39.15}{25.04}
\drawpath{47.15}{25.04}{59.15}{25.04}
\drawdot{41.15}{25.04}
\drawdot{45.15}{25.04}
\drawdot{43.15}{25.04}
\drawcenteredtext{19.15}{21.04}{\footnotesize 1}
\drawcenteredtext{27.15}{21.04}{\footnotesize 2}
\drawcenteredtext{51.15}{21.04}{\footnotesize $n-1$}
\drawcenteredtext{35.15}{21.04}{\footnotesize 3}
\drawcenteredtext{59.15}{21.04}{\footnotesize $n$}
\drawthickdot{73.08}{25.06}
\drawthickdot{81.08}{25.09}
\drawthickdot{89.08}{25.09}
\drawthickdot{105.08}{25.09}
\drawthickdot{113.08}{25.09}
\drawpath{73.08}{25.09}{89.08}{25.09}
\drawpath{89.08}{25.09}{93.08}{25.09}
\drawpath{101.08}{25.09}{113.08}{25.09}
\drawdot{95.08}{25.09}
\drawdot{99.08}{25.09}
\drawdot{97.08}{25.09}
\drawthickdot{73.08}{31.09}
\drawcenteredtext{81.08}{21.09}{\footnotesize 4}
\drawcenteredtext{105.08}{21.09}{\footnotesize $n-1$}
\drawcenteredtext{89.08}{21.09}{\footnotesize 5}
\drawcenteredtext{113.08}{21.09}{\footnotesize $n$}
\drawthickdot{73.08}{19.09}
\drawrighttext{71.08}{25.09}{\footnotesize 3}
\drawrighttext{71.08}{31.09}{\footnotesize 1}
\drawrighttext{71.08}{19.09}{\footnotesize 2}
\drawpath{73.08}{31.09}{73.08}{19.09}
\drawcenteredtext{40.21}{6.44}{\footnotesize (a)}
\drawthickdot{136.29}{37.27}
\drawthickdot{144.54}{34.37}
\drawthickdot{148.23}{25.0}
\drawthickdot{144.54}{16.09}
\drawthickdot{136.1}{12.93}
\drawthickdot{127.48}{16.25}
\drawthickdot{123.95}{24.84}
\drawthickdot{128.38}{34.47}
\drawdot{124.48}{29.59}
\drawdot{125.2}{30.75}
\drawdot{125.9}{31.89}
\drawpath{128.38}{34.47}{136.29}{37.29}
\drawpath{148.23}{25.0}{144.54}{16.09}
\drawpath{144.54}{16.09}{136.1}{12.93}
\drawpath{136.1}{12.93}{127.48}{16.25}
\drawpath{127.48}{16.25}{123.95}{24.84}
\drawcenteredtext{136.25}{40.00}{\footnotesize 1}
\drawcenteredtext{146.48}{36.14}{\footnotesize 2}
\drawcenteredtext{150.26}{25.10}{\footnotesize 3}
\drawcenteredtext{145.79}{14.24}{\footnotesize 4}
\drawcenteredtext{136.10}{10.47}{\footnotesize 5}
\drawcenteredtext{126.00}{14.14}{\footnotesize 6}
\drawcenteredtext{121.85}{25.10}{\footnotesize 7}
\drawcenteredtext{126.30}{36.14}{\footnotesize $n$}
\drawcenteredtext{90.13}{6.44}{\footnotesize (b)}
\drawcenteredtext{136.57}{6.44}{\footnotesize (c)}
\drawpath{136.38}{37.36}{144.54}{34.47}
\drawpath{144.54}{34.47}{148.32}{25.09}
\end{picture}
\vspace{-2em}
\caption{(a) the path $P_n$; (b) the T-caterpillar $T_n$; (c) the cycle $C_n$.\label{fig1}}
\end{figure}
\unitlength 1.50mm

All spanning rooted forests of $P_4$ and the spanning rooted forests in which vertex 1 belongs to a tree rooted at
vertex 2 are shown in Fig.~\ref{fig0}, where thick dots denote roots.
\unitlength 1.45mm  
\begin{figure}[htb] 
\begin{center}
\noindent
\begin{picture}(113.5,26)
\thicklines
\drawdot{6.0}{20.0}
\drawdot{6.0}{18.0}
\drawdot{6.0}{16.0}
\drawdot{6.0}{14.0}
\drawdot{10.0}{20.0}
\drawdot{10.0}{18.0}
\drawdot{10.0}{16.0}
\drawdot{10.0}{14.0}
\drawdot{14.0}{20.0}
\drawdot{14.0}{18.0}
\drawdot{14.0}{16.0}
\drawdot{14.0}{14.0}
\drawdot{18.0}{20.0}
\drawdot{18.0}{18.0}
\drawdot{18.0}{16.0}
\drawdot{18.0}{14.0}
\drawdot{22.0}{20.0}
\drawdot{22.0}{18.0}
\drawdot{22.0}{16.0}
\drawdot{22.0}{14.0}
\drawdot{26.0}{20.0}
\drawdot{26.0}{18.0}
\drawdot{26.0}{16.0}
\drawdot{26.0}{14.0}
\drawdot{30.0}{20.0}
\drawdot{30.0}{18.0}
\drawdot{30.0}{16.0}
\drawdot{30.0}{14.0}
\drawdot{34.0}{20.0}
\drawdot{34.0}{18.0}
\drawdot{34.0}{16.0}
\drawdot{34.0}{14.0}
\drawdot{34.0}{20.0}
\drawdot{34.0}{18.0}
\drawdot{34.0}{16.0}
\drawdot{34.0}{14.0}
\drawdot{38.0}{20.0}
\drawdot{38.0}{18.0}
\drawdot{38.0}{16.0}
\drawdot{38.0}{14.0}
\drawdot{42.0}{20.0}
\drawdot{42.0}{18.0}
\drawdot{42.0}{16.0}
\drawdot{42.0}{14.0}
\drawdot{46.0}{20.0}
\drawdot{46.0}{18.0}
\drawdot{46.0}{16.0}
\drawdot{46.0}{14.0}
\drawdot{50.0}{20.0}
\drawdot{50.0}{18.0}
\drawdot{50.0}{16.0}
\drawdot{50.0}{14.0}
\drawdot{54.0}{20.0}
\drawdot{54.0}{18.0}
\drawdot{54.0}{16.0}
\drawdot{54.0}{14.0}
\drawdot{58.0}{20.0}
\drawdot{58.0}{18.0}
\drawdot{58.0}{16.0}
\drawdot{58.0}{14.0}
\drawdot{62.0}{20.0}
\drawdot{62.0}{18.0}
\drawdot{62.0}{16.0}
\drawdot{62.0}{14.0}
\drawdot{66.0}{20.0}
\drawdot{66.0}{18.0}
\drawdot{66.0}{16.0}
\drawdot{66.0}{14.0}
\drawdot{70.0}{20.0}
\drawdot{70.0}{18.0}
\drawdot{70.0}{16.0}
\drawdot{70.0}{14.0}
\drawdot{74.0}{20.0}
\drawdot{74.0}{18.0}
\drawdot{74.0}{16.0}
\drawdot{74.0}{14.0}
\drawdot{78.0}{20.0}
\drawdot{78.0}{18.0}
\drawdot{78.0}{16.0}
\drawdot{78.0}{14.0}
\drawdot{82.0}{20.0}
\drawdot{82.0}{18.0}
\drawdot{82.0}{16.0}
\drawdot{82.0}{14.0}
\drawdot{86.0}{20.0}
\drawdot{86.0}{18.0}
\drawdot{86.0}{16.0}
\drawdot{86.0}{14.0}
\drawdot{94.0}{20.0}
\drawdot{94.0}{18.0}
\drawdot{94.0}{16.0}
\drawdot{94.0}{14.0}
\drawdot{98.0}{20.0}
\drawdot{98.0}{18.0}
\drawdot{98.0}{16.0}
\drawdot{98.0}{14.0}
\drawdot{102.0}{20.0}
\drawdot{102.0}{18.0}
\drawdot{102.0}{16.0}
\drawdot{102.0}{14.0}
\drawdot{106.0}{20.0}
\drawdot{106.0}{18.0}
\drawdot{106.0}{16.0}
\drawdot{106.0}{14.0}
\drawdot{110.0}{20.0}
\drawdot{110.0}{18.0}
\drawdot{110.0}{16.0}
\drawdot{110.0}{14.0}
\drawpath{10.0}{20.0}{10.0}{18.0}
\drawpath{14.0}{20.0}{14.0}{18.0}
\drawpath{18.0}{18.0}{18.0}{16.0}
\drawpath{22.0}{18.0}{22.0}{16.0}
\drawpath{26.0}{16.0}{26.0}{14.0}
\drawpath{30.0}{16.0}{30.0}{14.0}
\drawpath{34.0}{20.0}{34.0}{16.0}
\drawpath{38.0}{20.0}{38.0}{16.0}
\drawpath{42.0}{20.0}{42.0}{16.0}
\drawpath{42.0}{16.0}{42.0}{16.0}
\drawpath{46.0}{18.0}{46.0}{14.0}
\drawpath{50.0}{18.0}{50.0}{14.0}
\drawpath{54.0}{18.0}{54.0}{14.0}
\drawpath{58.0}{20.0}{58.0}{18.0}
\drawpath{58.0}{16.0}{58.0}{14.0}
\drawpath{62.0}{20.0}{62.0}{18.0}
\drawpath{62.0}{16.0}{62.0}{14.0}
\drawpath{66.0}{20.0}{66.0}{18.0}
\drawpath{66.0}{16.0}{66.0}{14.0}
\drawpath{70.0}{20.0}{70.0}{18.0}
\drawpath{70.0}{16.0}{70.0}{14.0}
\drawpath{74.0}{20.0}{74.0}{14.0}
\drawpath{78.0}{20.0}{78.0}{14.0}
\drawpath{82.0}{20.0}{82.0}{14.0}
\drawpath{86.0}{20.0}{86.0}{14.0}
\drawpath{94.0}{20.0}{94.0}{18.0}
\drawpath{98.0}{20.0}{98.0}{18.0}
\drawpath{98.0}{16.0}{98.0}{14.0}
\drawpath{102.0}{20.0}{102.0}{18.0}
\drawpath{102.0}{16.0}{102.0}{14.0}
\drawpath{106.0}{20.0}{106.0}{16.0}
\drawpath{110.0}{20.0}{110.0}{14.0}
\put(94.0,18.0){\circle*{0.8}}
\put(98.0,18.0){\circle*{0.8}}
\put(102.0,18.0){\circle*{0.8}}
\put(106.0,18.0){\circle*{0.8}}
\put(110.0,18.0){\circle*{0.8}}
\put(98.0,16.0){\circle*{0.8}}
\put(102.0,14.0){\circle*{0.8}}
\put(6.0,20.0){\circle*{0.8}}
\put(6.0,18.0){\circle*{0.8}}
\put(6.0,16.0){\circle*{0.8}}
\put(6.0,14.0){\circle*{0.8}}
\put(10.0,20.0){\circle*{0.8}}
\put(10.0,16.0){\circle*{0.8}}
\put(10.0,14.0){\circle*{0.8}}
\put(14.0,18.0){\circle*{0.8}}
\put(14.0,16.0){\circle*{0.8}}
\put(14.0,14.0){\circle*{0.8}}
\put(18.0,20.0){\circle*{0.8}}
\put(18.0,18.0){\circle*{0.8}}
\put(18.0,14.0){\circle*{0.8}}
\put(22.0,20.0){\circle*{0.8}}
\put(22.0,16.0){\circle*{0.8}}
\put(22.0,14.0){\circle*{0.8}}
\put(26.0,20.0){\circle*{0.8}}
\put(26.0,18.0){\circle*{0.8}}
\put(26.0,16.0){\circle*{0.8}}
\put(30.0,20.0){\circle*{0.8}}
\put(30.0,18.0){\circle*{0.8}}
\put(30.0,14.0){\circle*{0.8}}
\put(34.0,20.0){\circle*{0.8}}
\put(34.0,14.0){\circle*{0.8}}
\put(38.0,18.0){\circle*{0.8}}
\put(38.0,14.0){\circle*{0.8}}
\put(42.0,16.0){\circle*{0.8}}
\put(42.0,14.0){\circle*{0.8}}
\put(46.0,20.0){\circle*{0.8}}
\put(46.0,18.0){\circle*{0.8}}
\put(50.0,20.0){\circle*{0.8}}
\put(50.0,16.0){\circle*{0.8}}
\put(54.0,20.0){\circle*{0.8}}
\put(54.0,14.0){\circle*{0.8}}
\put(58.0,20.0){\circle*{0.8}}
\put(58.0,16.0){\circle*{0.8}}
\put(62.0,20.0){\circle*{0.8}}
\put(62.0,14.0){\circle*{0.8}}
\put(66.0,18.0){\circle*{0.8}}
\put(66.0,16.0){\circle*{0.8}}
\put(70.0,18.0){\circle*{0.8}}
\put(70.0,14.0){\circle*{0.8}}
\put(74.0,20.0){\circle*{0.8}}
\put(78.0,18.0){\circle*{0.8}}
\put(82.0,16.0){\circle*{0.8}}
\put(86.0,14.0){\circle*{0.8}}
\drawcenteredtext{3.5}{20.1}{{\scriptsize\bf 1}}
\drawcenteredtext{3.5}{18.1}{{\scriptsize\bf 2}}
\drawcenteredtext{3.5}{16.1}{{\scriptsize\bf 3}}
\drawcenteredtext{3.5}{14.1}{{\scriptsize\bf 4}}
\drawcenteredtext{91.5}{20.1}{{\scriptsize\bf 1}}
\drawcenteredtext{91.5}{18.1}{{\scriptsize\bf 2}}
\drawcenteredtext{91.5}{16.1}{{\scriptsize\bf 3}}
\drawcenteredtext{91.5}{14.1}{{\scriptsize\bf 4}}

\drawunderbrace{46.0}{11.0}{117.0}
\drawunderbrace{102.0}{11.0}{24.0}
\drawcenteredtext{46.0}{6.0}{\footnotesize $f=21$}
\drawcenteredtext{102.0}{6.0}{\footnotesize $f_{21}=5$}
\end{picture}
\end{center}
\vspace{-2.8em}
\caption{The spanning rooted forests in $P_4$ and the forests where 1 is in a tree rooted at~2.\label{fig0}}
\end{figure}
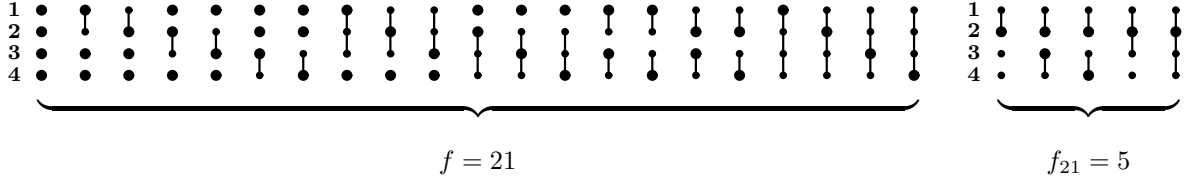
\unitlength 1.50mm

The matrix $F$ for $P_4$ is
$$
F(P_4)=\frac{(f_{ij})}{f}=
\frac{1}{21}\left[\begin{array}{rrrr}
13 & 5 & 2 & 1\\
5  &10 & 4 & 2\\
2  & 4 &10 & 5\\
1  & 2 & 5 &13
\end{array}\right].
$$

Let $T_n$ be the graph obtained from $P_n$ by replacing the edge $(1,2)$ with $(1,3)$:
$V(T_n)=\{1,2,\ldots,n\}$ and $E(T_n)=\{(1,3),(2,3),(3,4),\ldots,(n-1,n)\}$, see Fig.~\ref{fig1}(b). We call
$T_n$ a {\em T-caterpillar}.

Let $C_n$, $n\ge3$, be the {\em cycle\/} on $n$ vertices: $V(C_n)=\{1,2,\ldots,n\}$ and
$E(C_n)=\{(1,2),(2,3)\cdc(n-1,n),(n,1)\}$, Fig.~\ref{fig1}(c).

By $(\F_i)_{i=0,1,2,\ldots}=(0,1,1,2,3,5,\ldots)$ we denote the Fibonacci numbers. Sometimes, it is convenient to consider
the subsequences of Fibonacci numbers with odd and even subscripts separately:
\begin{eqnarray*}
\F'_i   &=&\F_{2i-1},\quad i=1,2,\ldots;\quad\\
\F''_{i}&=&\F_{2i},\quad\;\;\:\, i=0,1,2,\ldots.
\end{eqnarray*}

In Section~\ref{sec_nufor} we study the spanning rooted forests in paths, cycles, and T-caterpillars,
in Section~\ref{sec_proximity} the results are interpreted in terms of vertex-vertex proximity,
and Sections~\ref{sec_walk} and~\ref{sec_info} present interpretations of the doubly stochastic graph matrix
in terms of random walks and information dissemination, respectively.

\section{Spanning rooted forests in paths, cycles, and\\ {T-caterpillars}} 
\label{sec_nufor}

\begin{thm} 
\label{P_thm}
Let $G$ be a path$,$
$G=P_n$. Then $f=\F''_n$ and $f_{ij}=\F'_{\min(i,j)}\cdot\F'_{n+1-\max(i,j)}$
for all $i,j=1,\ldots,n$.
\end{thm}

\vspace{-0.3em} The number $f(G)$ of spanning rooted forests in any graph $G$ is equal to the number of
spanning trees in the graph $G^{+1}$, which is $G$ augmented by a ``hub'' vertex adjacent to every vertex
of~$G$. Indeed, a bijection between the spanning rooted forests of $G$ and spanning trees of $G^{+1}$ is
established by connecting every root of every spanning rooted forest to the ``hub'' by an edge. If $G=P_n$,
then $G^{+1}$ is a fan graph sometimes also called a ``terminated ladder.'' The fact that the number of
spanning trees in this $G^{+1}$ equals $\F''_n$ is familiar to electrical network theorists (\cite{Myers67},
cf.\,\cite{Morgan-Voyce59,Mowery61,Basin63MM}). 
Among others, it was obtained by Hilton~\cite{Hilton74}. Myers~\cite{Myers75} proved this using the notion of
weighted composition; in \cite{BoeschProdinger86} this fact was derived using Chebyshev polynomials. Our aim
is to give a direct combinatorial proof of Theorem~\ref{P_thm} in order to fully clarify the recurrence
structure of spanning rooted forests in a path. Here, a proof of $f(P_n)=\F''_n$ is integrated with a proof
of the expression for $f_{ij}$ given in Theorem~\ref{P_thm}.

For any graph $G$, $\FF(G)$ will denote the set of spanning rooted forests of~$G$.

\medskip\smallskip
\noindent{\bf Proof of Theorem~\ref{P_thm}.}
Let $\FF^m=\FF(P_m)$, let $f(m)=\abs{\FF^m}$, $m=1,2,\ldots.$ Then for every $k\ge1,$
\beq
\label{fn}
f(k+1)=\abs{\FF^{k+1}_{(1,2)}}+\abs{\FF^{k+1}_{(\overline{1,2})}},
\eeq
where
$\FF^{m}_{(1,2)}=\{\Ff\in\FF^{m}\,|\;(1,2)\in E(\Ff)\}$ and
$\FF^{m}_{(\overline{1,2})}=\FF^{m}\smallsetminus \FF^{m}_{(1,2)}$.

Let $\FF^m_*=\FF^{}_*(P_m)=\{\Ff\in\FF^m\,|\;1\text{ is a root in }\Ff\}$.
Then in (\ref{fn})
\beq
\label{f1to1}
\abs{\FF^{k+1}_{(\overline{1,2})}}=\abs{\FF^k}\;\text{  and  }\;
\abs{\FF^{k+1}_{(1,2)}}           =\abs{\FF_*^{k+1}}.
\eeq

Indeed, a bijection between $\FF^{k+1}_{(\overline{1,2})}$ and $\FF^k$ can be established by the
restriction of each $\Ff\in\FF^{k+1}_{(\overline{1,2})}$ to the vertex subset $\{2\cdc k\}$; a bijection between
$\FF^{k+1}_{(1,2)}$ and $\FF_*^{k+1}$ is established as follows: for every
$\Ff\in\FF^{k+1}_{(1,2)}$ obtain $\Ff'$ by
putting $\Ff'=\Ff$ if 1 is a root in~$\Ff$ and by
removing edge $(1,2)$ and marking vertex~1 as a root, otherwise.
Then $\Ff'\in\FF_*^{k+1}$ and this correspondence is one-to-one.

By (\ref{fn}) and (\ref{f1to1}),
\beq
\label{fn1}
f(k+1)=f(k)+f^*(k+1)\quad k=1,2,\ldots,
\eeq
where $f^*(m)=\abs{\FF^m_*}$.
Let $\FF^{m}_{(*1,2)}=\FF^{m}_{*}\cap\FF^{m}_{(1,2)}$, $m=1,2,\ldots.$ Using (\ref{f1to1}) we obtain
\beq
\label{f*n}
f^*(k+1)=\abs{\FF^{k+1}_{(\overline{1,2})}}+\abs{\FF^{k+1}_{(*1,2)}}=\abs{\FF^k}+\abs{\FF^k_*}=f(k)+f^*(k),
\quad k=1,2,\ldots.
\eeq

Here, a bijection between $\FF^{k+1}_{(*1,2)}$ and $\FF^k_*$ is established by coalescing 
vertex~2 with the root~1 and collapsing edge $(1,2)$.

Observe now that $f(1)=1=\F''_1$ and $f^*(1)=1=\F'_1$. By (\ref{fn1}) and (\ref{f*n}), $f(k)$ and $f^*(k)$ satisfy
the same recurrence relations as $\F''_k$ and $\F'_k$ do, respectively. Therefore
\beq
\label{equiv}
f(k)=\F''_k \;\:\text{ and }\:\; f^*(k)=\F'_k,\quad k=1,2,\ldots.
\eeq

Thus, $f=f(n)=\F''_n$. To find $f_{ij}$, $i,j=1\cdc n$, observe that $f_{ij}$ counts the spanning rooted
forests that contain the $i$--$j$ path rooted at~$i$. To obtain a spanning rooted forest, this path can be
completed on the subset of vertices $\{1\cdc\min(i,j)\}$ in $f^*(\min(i,j))$ ways and on the subset of
vertices $\{\max(i,j)\cdc n\}$ in $f^*(n+1-\max(i,j))$ ways. Since the ways of these types are all
compatible, (\ref{equiv}) implies that $f_{ij}=\F'_{\min(i,j)}\cdot\F'_{n+1-\max(i,j)}$. \qed
\bigskip

Theorem~\ref{P_thm} as well as Theorem~\ref{C_thm} below can also be proved algebraically by means of the matrix
forest theorem (Eqs.~(\ref{mft1}) and~(\ref{mft2})). We present combinatorial proofs here, since they are a
bit more illuminating. However, Theorem~\ref{T_thm} below is proved algebraically.

\begin{thm}
\label{C_thm} Let $G$ be a cycle$,$ $G=C_n$ with $n\ge3$. Then $f=\F'_n+\F'_{n+1}-2$ and\linebreak
$f_{ij}=\F''_{\abss{j-i}}+\F''_{n-\abss{j-i}},\;i,j=1\cdc n$.
\end{thm}

For $G=C_n$, the augmented graph $G^{+1}$ mentioned above is the wheel on $n+1$ vertices. The fact that the
number of spanning trees in the wheel is $\F'_n+\F'_{n+1}-2$ is due to Sedl\'a\v{c}ek \cite{Sedlacek70} and
Myers~\cite{Myers71}. Myers~\cite{Myers75} proved this using identities involving weighted compositions; the
proof by Benjamin and Yerger \cite{BenjaminYerger06} is based on counting imperfect matchings. A useful tool
for solving such problems is Chebyshev polynomials, see~\cite{Mowery61,BoeschProdinger86,ZhangYongGolin05}. Our proof
of the identity $f(C_n)=\F'_n+\F'_{n+1}-2$ presented here for completeness is based on relations between
forests found before. The proof of Theorem\!~\ref{C_thm} relies on the following lemma.

\begin{lemma}
\label{l_patpat}
For $n\ge2,$ let $\FF^{n}_{**}=\{\Ff\in\FF(P_{n})\,|\,1\text{ and }\, n \text{ are roots in } \Ff\}$.
Then $\abs{\FF^{n}_{**}}=\F''_{n-1}$.
\end{lemma}

\noindent{\bf Proof.}
A bijection between $\FF^{n}_{**}$ and the set $\FF(P_{n-1})$ of spanning rooted forests in $P_{n-1}$ can be
established as follows. For every $\Ff\in\FF(P_{n-1})$ define $\Ff'$ as the spanning subgraph of $P_n$ whose
roots satisfy two conditions:

(1) vertex $i$ is a root in $\Ff'$ iff [$i=1$ or $i=n$ or $(i-1,i)\not\in E(\Ff)$];

(2) $(i,i+1)\in E(\Ff')$ iff $i$ is not a root in $\Ff$.
\smallskip

In this case, $\Ff'$ is a spanning rooted forest of $P_n$. Indeed, if one assumes that some tree in $\Ff'$
has no root or has more than one root, then this would imply the presence of a tree with more than one root
or no root in~$\Ff$, respectively. Furthermore, every $\Ff'\in\FF^{n}_{**}$ has a pre-image in
$\FF(P_{n-1})$, and this correspondence is by definition one-to-one.

Consequently, by (\ref{equiv}), $\abs{\FF^{n}_{**}}=\abss{\FF(P_{n-1})}=\F''_{n-1}$.
\qed
\bigskip

\noindent{\bf Proof of Theorem~\ref{C_thm}.}
Let $\FF^m_{ij}$ be the set of spanning rooted forests of $C_m$ in which $j$ belongs to a tree rooted at~$i$.
Observe that
\beq
\label{spl_cyc}
f_{ij}=\abss{\FF^n_{i\frown j}}+
       \abss{\FF^n_{i\smile j}},
\eeq
where
$\FF^n_{i\frown j}=\{\Ff\in\FF^n_{ij}\,|\,(n,1)\not\in E(\Ff)\}$ and
$\FF^n_{i\smile j}=\{\Ff\in\FF^n_{ij}\,|\,(n,1)    \in E(\Ff)\}$.

We now show that $\abss{\FF^n_{i\smile j}}=\F''_{\abss{j-i}}$.
Every forest in $\FF^n_{i\smile j}$ contains the path $P_{n+1-\abss{j-i}}$ formed by the vertices in the
sequence $\left(\max(i,j)\cdc n,1\cdc\min(i,j)\right)$ and the edges between the neighboring elements in this sequence.
The ways of completing this path to obtain a spanning rooted forest in $C_n$ can be put into a one-to-one correspondence
with the elements of the set $\FF^{\abss{j-i}+1}_{**}$ defined in Lemma~\ref{l_patpat}. Indeed, linking each
$\tilde \Ff\in\FF^{\abss{j-i}+1}_{**}$ with $P_{n+1-\abss{j-i}}$ by replacing the vertices 1 and $\abss{j-i}+1$ of
$\tilde \Ff$ with vertices $i$ and $j$ of $P_{n+1-\abss{j-i}}$, respectively, produces a spanning rooted forest of $C_n$,
and every spanning rooted forest of $C_n$ can be uniquely obtained in this manner. Thus by Lemma~\ref{l_patpat},
$\abss{\FF^n_{i\smile j}}=\abs{\FF^{    \abss{j-i}+1}_{**}}=\F''_{  \abss{j-i}}$. Similarly,
$\abss{\FF^n_{i\frown j}}=\abs{\FF^{n+1-\abss{j-i}  }_{**}}=\F''_{n-\abss{j-i}}$.
Therefore by (\ref{spl_cyc}) it follows that $f_{ij}=\F''_{\abss{j-i}}+\F''_{n-\abss{j-i}}$.

Let $\FF^m=\FF(C_m)$, $m\ge3$. Then
\beq
\label{Cfdeco}
f(C_n)=\abs{\FF_{(\on)}}+\abs{\FF_{(1,n)*}}+\abs{\FF_{*(1,n)}},
\eeq
where
\vspace{-1.0em}
\beqq
\FF_{(\on)}&\!=\!&\{\Ff\in\FF^m\,|\,(1,n)\not\in E(\Ff)\},\nonumber\\ \nonumber
\FF_{          (1,n)*}&\!=\!&\{\Ff\in\FF^m\,|\,(1,n)\in E(\Ff)\text{ and the path joining 1 with the root contains }n\},\\
\FF_{         *(1,n) }&\!=\!&\{\Ff\in\FF^m\,|\,(1,n)\in E(\Ff)\text{ and the path joining }n
                                                              \text{ with the root contains 1}\}. \nonumber
\eeqq

Obviously, $\abs{\FF_{(\on)}}=\abss{\FF(P_n)}$, so, by Theorem~\ref{P_thm}, $\abs{\FF_{(\on)}}=\F''_n$.
Consider any $\Ff\in\FF_{(1,n)*}$. Removing $(1,n)$ from $E(\Ff)$ and marking 1 as a root produces a forest
$\Ff'\in\FF_*(P_n)$, where $\FF^{}_*(P_n)=\{\Ff\in\FF(P_n)\,|\;1\text{ is a root in }\Ff\}$ was defined in
the proof of Theorem~\ref{P_thm}. All elements of $\FF_*(P_n)$ can be obtained in this way, except for the
whole path $P_n$ rooted at~1. That is why $\abs{\FF_{(1,n)*}}=\abss{\FF_*(P_n)}-1$ and, by~(\ref{equiv}),
$\abs{\FF_{(1,n)*}}=\F'_n-1$. Similarly, $\abs{\FF_{*(1,n)}}=\F'_n-1$. Substituting this in (\ref{Cfdeco}) provides
\[
f(C_n)=\F''_n+2(\F'_n-1)=\F'_n+\F'_{n+1}-2. \eqno{\square}
\]

\bigskip
Now recall that
\beq
\nonumber
\La_i=\F_{i-1}+\F_{i+1},
\eeq
where $\F_{-1}=1$,
are the Lucas numbers: $(\La_i)_{i=0,1,2,\ldots}=(2,1,3,4,7,11,18,29,47,\ldots)$, see, e.g.,~\cite{Koshy01}.
The Lucas numbers satisfy the same recurrence as the Fibonacci numbers do:
\[
\La_i+\La_{i+1}=\La_{i+2}, \quad i=0,1,2,\ldots,
\]
but some other properties of the Lucas numbers are even more elegant than those of the Fibonacci numbers.

By Theorem~\ref{C_thm}, $f(C_n)=\La_{2n}-2$. The numbers of forests in a cycle can also be expressed via
smaller Fibonacci and Lucas numbers, viz. Corollary~\ref{C_cor1} holds.\footnote{A knot theory interpretation
of the squareness of $\La_{2n}-2$ when $n$ is odd can be found in \cite{Stoimenow-wheel00}.}

\begin{corol}
\label{C_cor1} Let $G$ be a cycle$,$ $G=C_n$ with $n\ge3$. Then
\[
f=
\begin{cases}
\La_n^2, &n=2k-1\\
5\F_n^2, &n=2k
\end{cases};
\quad\;\:
f_{ij}=
\begin{cases}
 \F_t\La_n, &n=2k-1\\
\La_t \F_n, &n=2k
\end{cases},
\quad i,j=1\cdc n,
\]
where $t=\bigl|n-2\abss{j-i}\bigr|$.
\end{corol}

Corollary~\ref{C_cor1} is derived from Theorem~\ref{C_thm} by means of classical identities involving
Fibonacci and Lucas numbers.
It provides a simple expression for the entries of the doubly stochastic graph matrix
$F=\frac{(f_{ij})}{f}$ of $C_n$:
\begin{corol}
\label{C_cor2}
The entries of the doubly stochastic matrix $F=\frac{(f_{ij})}{f}$ of\/ $C_n$\/ $(n\ge3)$ are$:$
\beq
\label{FCn}
\frac{f_{ij}}{f}=
\begin{cases}
 \F_{t}/\La_n, &n=2k-1\\
\La_{t}/5\F_n, &n=2k
\end{cases},
\quad i,j=1\cdc n,
\eeq
where $t=\bigl|n-2\abss{j-i}\bigr|$.
\end{corol}

In the expression (\ref{FCn}), for every row of $F$, the numerators make up a segment of a fixed symmetric two-sided
sequence:
for odd  $n$ this sequence is $(\ldots,34,13,5,2,1,1,2,5,13,34,\ldots)$;
for even $n$            it is $(\ldots,47,18,7,3,2,3,7,18,47,\ldots).$
Thereby the ratio of two corresponding elements of $F$ is the same for all $n$ of the same parity.

\medskip
Regarding the T-caterpillars, we are mainly interested in the total number $f$ of spanning rooted
forests and the diagonal entries $f_{33}$ and $f_{nn}$ of the matrix of spanning rooted forests.
\begin{thm}
\label{T_thm}
\,Let $G$ be a T-caterpillar$,$ $G\!=\!T_n$.\! Then $f\!=\!4\F'_{n-1},$ $f_{33}\!=\!4\F'_{n-2},$ and $f_{nn}\!=\!4\F''_{n-2}$.
\end{thm}

\noindent{\bf Proof.}
Observe that $\F''_0=0$, $\F''_1=1$, and for $i=1,2,\ldots,$
\beqq
\label{recuF''}
\nonumber
\F''_{i+1}&=&\F_{2i+2}=\F_{2i}+\F_{2i+1}=
            \F_{2i}+\F_{2i}+\F_{2i-1}\\
          &=&
            2\F_{2i}+\F_{2i}-\F_{2i-2}=
            3\F''_{i}-\F''_{i-1}.
\eeqq
Similarly, $\F'_1=1$ and for $i=1,2,\ldots,$
\beq
\label{recuF'}
\F'_{i+1}=
\begin{cases}
3\F'_{i}-\F'_{i-1},&i>1\\
2\F'_{i},          &i=1
\end{cases}.
\eeq

For a T-caterpillar,
\beq
\label{Tmatr}
I+L=
\left[\begin{array}{rrrrrrrrr}
  2 & 0 & -1& 0 & 0 &\cdots&  0 & 0 & 0\\
  0 & 2 & -1& 0 & 0 &\cdots&  0 & 0 & 0\\
  -1& -1& 4 & -1& 0 &\cdots&  0 & 0 & 0\\
  0 & 0 & -1& 3 & -1&\cdots&  0 & 0 & 0\\
  0 & 0 & 0 & -1& 3 &\cdots&  0 & 0 & 0\\
  \hdotsfor[2]{9}                      \\
  0 & 0 & 0 & 0 & 0 &\cdots&  3 & -1& 0\\
  0 & 0 & 0 & 0 & 0 &\cdots&  -1& 3 &-1\\
  0 & 0 & 0 & 0 & 0 &\cdots&  0 & -1& 2\\
\end{array}\right].
\eeq
Equations $F(I+L)=I$ (see (\ref{mft1})) and (\ref{Tmatr}) imply
\begin{eqnarray}
2f_{n1}-f_{n3}=0,               &           \nonumber\\
2f_{n2}-f_{n3}=0,               &           \nonumber\\
-f_{n1}-f_{n2}+4f_{n3}-f_{n4}=0,&           \nonumber\\
-f_{n3}+3f_{n4}-f_{n5}=0,       &           \label{eqar3}\\
........................................&   \nonumber\\
-f_{n,n-2}+3f_{n,n-1}-f_{nn}=0, &           \nonumber\\
-f_{n,n-1}+2f_{nn}          =f. &           \nonumber
\end{eqnarray}

Note that $f_{n1}=f_{n2}=2=2\F''_1$, since $T_n$ has exactly two spanning rooted forests where vertex
$1$ belongs to a tree rooted at~$n$, one of them being $T_n$ with root $n$, the other $T_n$ with edge
$(2,3)$ deleted and roots $n$ and~$2$. Then, by (\ref{eqar3}), $f_{n3}=2f_{n1}=4=4\F''_{1}$. Consequently,
using (\ref{eqar3}), (\ref{recuF''}) and induction, we obtain that
\beqq
&&f_{n4}=\,4f_{n3}-f_{n1}-f_{n2}=12=4(3\F''_{1}-\F''_{0})=4\F''_{2},\nonumber\\
&&f_{n5}=\,3f_{n4}-f_{n3}=4(3\F''_{2}-\F''_{1})=4\F''_{3},\nonumber\\
&&........................................................................\label{TnInfer1}  \\
&&f_{nn}=\,3f_{n,n-1}-f_{n,n-2}=4(3\F''_{n-3}-\F''_{n-4})=4\F''_{n-2}.\nonumber
\eeqq
From the last equation of (\ref{eqar3}),
$f=2f_{nn}-f_{n,n-1}=4(2\F''_{n-2}-\F''_{n-3})=4(\F''_{n-2}+\F'_{n-2})=4\F'_{n-1}$.
It remains to show that $f_{33}=4\F'_{n-2}$. Eqs.~(\ref{mft1}) and (\ref{Tmatr}) imply
\begin{eqnarray}
-f_{33}+3f_{34}-f_{35}=0,               &   \nonumber\\
-f_{34}+3f_{35}-f_{36}=0,               &   \nonumber\\
........................................&   \label{eqar4}\\
-f_{3,n-2}+3f_{3,n-1}-f_{3n}=0,         &   \nonumber\\
-f_{3,n-1}+2f_{3n}=0.                   &   \nonumber
\end{eqnarray}
Since $f_{3n}=f_{n3}=4=4\F'_{1}$, from (\ref{eqar4}) and (\ref{recuF'}) we have that $f_{3,n-1}=4\F'_{2}$ and, by
induction, $f_{33}=4\F'_{n-2}$.
\qed

\section{``Golden introverts'' and ``golden extroverts''}
\label{sec_proximity}

In \cite{CheSha95,CheSha97,CheSha98} (see also \cite{Merris97,Merris98,CheSha01})
$F=\frac{(f_{ij})_{n\times n}}{f}$ was studied as a  matrix of vertex-vertex proximity. For every graph $G$,
$F$~is a positive definite
doubly stochastic matrix and $\frac{f_{ij}}{f}$ measures the relative strength of connections between vertices $i$ and
$j$\/ in~$G$. This proximity measure was referred to as the {\em relative connectivity via forests}. For some additional
applications of $F$ we refer to~\cite{GolenderDrboglav81,FoussPirotte07simi}.

It turns out \cite[Theorem~2]{Merris98} that for every pair of vertices $i$ and $j$ such that $j\ne i$,
$f_{ij}\le f_{ii}/2$; $\frac{f_{ii}}{f}$ can be considered as a measure of self-connectivity of vertex~$i$.
By \cite[Corollary~7]{Merris97}, $\frac{f_{ii}}{f}\ge (1+d_i)^{-1}$, where $d_i$ is the degree of vertex~$i$.

A vertex $i$ can be called ``an introvert'' if $\frac{f_{ii}}{f}>0.5$ (or, equivalently, $f_{ii}>\sum_{j\ne
i}f_{ij}$) and ``an extrovert'' if $\frac{f_{ii}}{f}<0.5$ (equivalently, $f_{ii}<\sum_{j\ne i}f_{ij}$). The
complete graph on three vertices provides an example of the boundary case where $\frac{f_{ii}}{f}=0.5$
and $f_{ii}=\sum_{j\ne i}f_{ij}$ for every vertex~$i$.

\begin{prop}
\label{goldenPro}
Let $\phi$ be the golden ratio$,$ $\phi=\frac{\sqrt{5}+1}{2}$. Then\newline
$(i)$  For the          paths $P_n,$ $\liml_{n\to\infty}\frac{f_{11}}{f}=  \phi^{-1}\,;$\newline
$(ii)$ For the T-caterpillars $T_n,$ $\liml_{n\to\infty}\frac{f_{nn}}{f}=  \phi^{-1}$ and\/
                                     $\liml_{n\to\infty}\frac{f_{33}}{f}=1-\phi^{-1}$.
\end{prop}

\noindent{\bf Proof.}
By Theorem~\ref{P_thm} and Binet's Fibonacci number formula, for the paths $P_n$,
                        $\liml_{n\to\infty}({f_{11}}/{f})=
                         \liml_{n\to\infty}({\F'_n}/      {\F''_n})    =\phi^{-1}$.
By Theorem~\ref{T_thm}, for the T-caterpillars $T_n$,
                        $\liml_{n\to\infty}({f_{nn}}/{f})=
                         \liml_{n\to\infty}({4\F''_{n-2}}/{4\F'_{n-1}})=\phi^{-1}$ and
                        $\liml_{n\to\infty}({f_{33}}/{f})=
                         \liml_{n\to\infty}({4\F'_{n-2}}/ {4\F'_{n-1}})=\phi^{-2}=1-\phi^{-1}$.
\qed
\medskip

\begin{corol}
\label{goldenProC}
Let $\phi$  be the golden ratio$,$ $\phi=\frac{\sqrt{5}+1}{2}$. Then\newline
$(i)$  For the          paths   $P_n,$  $\liml_{n\to\infty}\frac{f_{11}}{\text{\footnotesize$\mathstrut$}\suml_{i\ne 1}f_{1i}}=\phi\,;$\newline
$(ii)$ For the T-caterpillars\/ $T_n,$  $\liml_{n\to\infty}\frac{f_{nn}}{\text{\footnotesize$\mathstrut$}\suml_{i\ne n}f_{ni}}=\phi\,$ and\/\,
                                    $\liml_{n\to\infty}\frac{\suml_{i\ne 3}f_{3i}}{f_{33}}=\phi$.
\end{corol}

Corollary~\ref{goldenProC} follows from Proposition~\ref{goldenPro} and the fact that $F$ is stochastic.

It can be shown that (ii) of Proposition~\ref{goldenPro} and (ii) of Corollary~\ref{goldenProC} remain true for the
graphs resulting from T-caterpillars by the addition of edge $(1,2)$.

In accordance with Corollary~\ref{goldenProC}, as ${n\to\infty}$, vertices 1 and $n$ in a path and vertex $n$
in a T-caterpillar tend to be ``golden introverts'' (named after the golden ratio), whereas vertex 3 in a
T-caterpillar tends to become a ``golden extrovert.'' This provides a kind of  sociological interpretation of
Corollary~\ref{goldenProC}.

\section{A random walk interpretation of the doubly\\ stochastic graph matrix}
\label{sec_walk}

To better comprehend what exactly the results of the previous section mean, consider a random walk interpretation of
the doubly stochastic graph matrix.

For a graph $G$, consider any Markov chain whose states are the vertices of $G$, $\{1,2\cdc n\}$, and the
probabilities of all $i\!\to\!j$ transitions with $i\ne j$ are proportional\footnote{There are two popular
methods of attaching a Markov chain to a graph. The first one is based on~(\ref{graphMarkov}); for any undirected graph
it provides a symmetric transition matrix with, in general, nonzero diagonal. The second method assumes that
$p_{ij}=a_{ij}/\sum_{k=1}^na_{ik}$. For an undirected graph without loops it generally provides a nonsymmetric transition
matrix with zero diagonal.}
to the corresponding elements of the adjacency matrix of $G$:
\beq
\label{graphMarkov}
p_{ij}=\e a_{ij}, \quad i,j=1\cdc n,\;\;i\ne j.
\eeq
Then the diagonal elements of the transition matrix $P=(p_{ij})$ are determined as follows:
\beq
\label{dia_trans}
p_{ii}=1-\sum_{j\ne i}\e a_{ij},\quad i=1\cdc n
\eeq
and, in a matrix form,
\beq
\nonumber
P=I-\e L(G),
\eeq
where $L(G)$ is the Laplacian matrix of~$G$.

The maximum value of $\e$ that guarantees correctness, i.e., the nonnegativity of the diagonal entries
(\ref{dia_trans}) for all simple graphs on $n$ vertices, is obviously $\e=(n-1)^{-1}$. On the other hand,
$\e=(n-1)^{-1}$ is the only correct $\e$ that allows the self-transition  probabilities $p_{ii}$ to be zero.
Therefore it makes sense to consider this value of $\e$ and the Markov chain with transition matrix
\beq
\label{transM1}
P=I-(n-1)^{-1}L(G)
\eeq
in more detail.

For this chain, let us examine random walks with a random number of steps.
Namely, consider 
a sequence of independent Bernoulli trials indexed by $0,1,2,\ldots$ with a certain success probability~$q$.
Suppose that the number of steps in a random walk equals the trial number of the first success.
Then the number of steps, $K$, is a geometrically distributed random variable:
\[
\Pr\{K=k\}=q(1-q)^k,\quad k=0,1,2,\ldots.
\]

Suppose that $q={1}/{n}$. For this value, the expected number of steps is $n-1$, which is the number
of edges in every spanning tree of~$G$. Then
\beq
\label{geom_dis}
\Pr\{K=k\}=\frac{1}{n}\left(1-\frac{1}{n}\right)^k,\quad k=0,1,2,\ldots.
\eeq
Let $Q=(q_{ij})$ be the matrix with entries
\beq
\label{Qtran}
q_{ij}=\Pr\{X_K=j\mid X_0=i\},\quad i,j=1\cdc n,
\eeq
where $X_k$ is the state of the Markov chain under consideration at step $k$, i.e., $Q$ is the transition matrix
of the overall random walk with a random number of steps~$K$.

\begin{thm}
\label{prop_walk}
For a graph $G$ on $n$ vertices and the corresponding Markov chain whose transition matrix is\/ {\rm(\ref{transM1}),}
let $Q$ be the transition matrix\/ {\rm(\ref{Qtran})} of the overall random walk whose number of steps is geometrically
distributed with parameter~$1/n$. Then $Q=F,$ where $F=\frac{(f_{ij})_{n\!\times\!n}}{f}$ is the
doubly stochastic matrix of\/~$G$. 
\end{thm}

\noindent{\bf Proof.} Since the spectral radius of $P$ is 1, for every $q$ such that $0<s<1$
\[
\suml_{k=0}^{\infty}\left(sP\right)^{k}=\left(I-sP\right)^{-1}
\]
holds. Consequently, using the formula of total probability, (\ref{geom_dis}), (\ref{transM1})
and the matrix forest theorem~(\ref{mft1}) we obtain
\begin{eqnarray*}
Q&=& \sum_{k=0}^{\infty}\Pr\{K=k\}\,P^k\,
  =\,\sum_{k=0}^{\infty}\frac{1}{n}\left(1-\frac{1}{n}\right)^kP^k\\
 &=&                    \frac{1}{n}\left(I-\left(1-\frac{1}{n}\right)P\right)^{-1}
  =(I+L)^{-1}
  =F.
\end{eqnarray*}

\vspace{-2.45em}\qed

\bigskip\medskip
By virtue of Theorem~\ref{prop_walk}, if a ``golden extrovert'' walks randomly in accordance with
the above model, she eventually finds herself on a visit $\phi\,$ times more often than at home, whereas for a
``golden introvert'' the situation is opposite.

\def\baselinestretch{0.81}
\section{A concluding note: a communicative interpretation of the doubly stochastic graph matrix}
\label{sec_info}

In closing, let us mention an interpretation of the doubly stochastic graph matrix in terms of information dissemination.
Suppose that a sequence of information units (or ideas) are transmitted through a graph~$G$.
A {\em plan\/} of information transmission is a rooted forest $\Ff\in\FF(G)$:
every information unit (idea) is initially injected into the roots of $\Ff$; after that it comes to the other vertices
along the edges of~$\Ff$. Suppose that every time a possible plan is
chosen at random: the probability of every choice is $\abss{\FF(G)}^{-1}=f.$ Then
$\frac{(f_{ij})_{n\!\times\!n}}{f}$
is the probability that an information unit arrives at $j$ {\em from root}~$i$. As a result, for a ``golden introvert''
the expected proportion of ``her own'' (injected straight into her mind) ideas to adopted ideas is $\phi$, whereas for a
``golden extrovert'' the proportion is inverse.



\end{document}